\documentclass[11pt]{article}
\usepackage{graphics}
\usepackage{graphicx}

\RequirePackage[OT1]{fontenc}
\RequirePackage{amsmath,natbib}
\RequirePackage[colorlinks]{hyperref}
\RequirePackage{hypernat}

\usepackage[psamsfonts]{amsfonts}
\usepackage{pdflscape}

\newtheorem{Prop}{Proposition}[section]

\newcommand{\cqfd}{\hfill $\square$}

\newcommand{\R}{\mathbb R}

\newcommand{\n}{^{(n)}}

\newcommand{\Xb}{\mathbf{X}}
\newcommand{\Sb}{\mathbf{S}}

\newcommand{\zerob}{\ensuremath{\mathbf{0}}}
\newcommand{\xb}{\ensuremath{\mathbf{x}}}

\newcommand{\tb}{\mathbf{t}}

\newcommand{\thetab}{{\pmb \theta}}

\newcommand{\varthetab}{{\pmb \vartheta}}

\newcommand{\Umb}{{\pmb \Upsilon}}

\newcommand{\Deltab}{{\pmb \Delta}}
\newcommand{\taub}{{\pmb \tau}}
\newcommand{\nub}{{\pmb \nu}}

\newcommand{\etab}{{\pmb \eta}}

\newcommand{\Gamb}{{\pmb \Gamma}}

\newcommand{\pr}{^{\prime}}

\newcommand{\ny}{n\rightarrow\infty}

\begin{document}

\title{Local powers of optimal one- and multi-sample tests for the concentration of Fisher-von Mises-Langevin distributions 
}

\author{Christophe Ley and Thomas Verdebout}


\maketitle

\begin{abstract}

One-sample and multi-sample tests on the concentration parameter of Fisher-von Mises-Langevin (FvML) distributions have been well studied in the literature. However, only very little is known about their behavior under local alternatives, which is due to complications inherent to the curved nature of the parameter space. The aim of the present paper therefore consists in filling that gap by having recourse to the Le Cam methodology, which has been adapted from the linear to the spherical setup in Ley \emph{et al.}~(2013a). We obtain explicit expressions of the powers for the most efficient one- and multi-sample tests; these tests are those considered in Watamori and Jupp~(2005). As a nice by-product, we are also able to write down the powers (against local FvML alternatives) of the celebrated Rayleigh~(1919) test of uniformity. A Monte Carlo simulation study confirms our theoretical findings and shows the finite-sample behavior of the above-mentioned procedures.

\end{abstract}


\section{Introduction}
\label{intro}

The field of directional (circular in dimension $k=2$, spherical in higher dimensions) statistics has become increasingly popular over the past decades, stimulated by the pioneering and seminal paper Fisher~(1953). More recent cornerstone references are the monographs Fisher \emph{et al.}~(1987) and Mardia and Jupp~(2000). This domain, which is particularly suited for modeling and explaining phenomena arising in earth sciences, meteorology, the study of animal behavior, astronomy or neurosciences, to cite but these areas, is concerned with observations conceived as realizations of random vectors $\Xb$  taking values on the surface of the unit hypersphere $\mathcal{S}^{k-1}:= \{{\bf v} \in \R^k\,|\, \| {\bf v} \|=1\}$.

By far the most popular and most used directional distribution is the \emph{Fisher-von Mises-Langevin} (FvML) distribution
(named, according to Watson~1983, after von Mises~1918 for $k=2$, Fisher~1953 for $k=3$,
and Langevin~1905 for general $k$), whose  density is of the form (with respect to the usual surface area measure on spheres)
\begin{equation}\label{density}
f_{\kappa,\thetab}({\bf x})= c_{k,\kappa}  \exp(\kappa {\bf x}\pr \thetab),\quad\xb\in\mathcal{S}^{k-1},\end{equation}
where $\kappa>0$ is a concentration parameter,
$\thetab\in\mathcal{S}^{k-1}$ a (spherical or directional) location parameter and where the
normalizing constant $c_{k,\kappa}$ is given by 
$$
c_{k,\kappa}=\left(\frac{\kappa}{2}\right)^{k/2-1}\frac{1}{\Gamma(k/2)I_{k/2-1}(\kappa)},
$$
with $\Gamma(k/2)$ the Gamma function evaluated at $k/2$ and $I_{k/2-1}(\kappa)$ the modified Bessel function of the first kind and of order $k/2-1$.  The FvML distribution is considered as the directional
analogue of the (linear) Gaussian distribution for purposes of mathematical
statistics (see Schaeben~1992 for a discussion on directional analogues of the
Gaussian distribution). This analogy is mainly due to the fact that
the FvML distribution can be characterized by the empirical spherical
mean $\hat{\thetab}_{\rm
  Mean}:=\sum_{i=1}^n\Xb_i/||\sum_{i=1}^n\Xb_i||$,
$\Xb_1,\ldots,\Xb_n\in\mathcal{S}^{k-1}$, as the Maximum Likelihood
Estimator (MLE) of its (spherical) location parameter, similarly as the
Gaussian distribution can be characterized by the empirical mean
$n^{-1}\sum_{i=1}^n\Xb_i$, $\Xb_1,\ldots,\Xb_n\in\R^k$, as the MLE of
its classical (linear) location parameter, a famous result due to Gauss. We refer to Duerinckx and Ley~(2013) for a formal proof of this fact and for more details on directional MLE characterizations\footnote{It is interesting in this context to note that
  Gauss, in his manuscript ``Theoria motus corporum coelestium in
  sectionibus conicis solem ambientium'' of 1809, has defined the
  famous distribution named after him by searching for the probability law
  for which the sample mean is always the MLE of the location
  parameter, and that von Mises, in 1918, aiming at constructing a
  circular analogue of the Gaussian distribution, started precisely
  from this classical MLE characterization.}.
  
Due to its prominent role, the FvML distribution has received a lot of attention in the literature, and inferential procedures involving its concentration and location parameters have been extensively studied in the literature (see for instance Sections 10.4-10.6 in Mardia and Jupp~2000). In the present paper, the parameter of interest is the concentration parameter $\kappa$ which regulates the probability mass in the vicinity of the modal direction $\thetab$. Besides the tests described in Mardia and Jupp~(2000), hypothesis testing procedures dealing specifically with the concentration parameter can mainly be found in Stephens (1969), Larsen \emph{et al}. (2002) and Watamori and Jupp (2005). Due to their efficiency properties, the proposed procedures are either likelihood ratio  (in its basic and improved versions) or score tests. Even if the asymptotic theory of such tests has been well studied in the above-cited papers, little is known about their asymptotic behavior and power under local alternatives. This absence of result can certainly be explained by the curved nature of the parameter space $\R_0^+\times\mathcal{S}^{k-1}$. 

In this paper, our aim is therefore to fill this gap by providing explicit expressions of the powers of the most efficient tests for both the one-sample problem ($\mathcal{H}_0^{\kappa_0}: \kappa= \kappa_0$ for some fixed $\kappa_0 >0$ versus $\mathcal{H}_1^{\kappa_0}: \kappa\neq \kappa_0$) and the multi-sample problem ($\mathcal{H}_0^{\rm Hom}: \kappa_1= \ldots= \kappa_m$ for $m\geq2$ and $\kappa_1,\ldots,\kappa_m>0$ versus $\mathcal{H}_1^{\rm Hom}: \exists 1\leq i\neq j\leq m\,\, \kappa_i\neq \kappa_j$). We achieve this goal by combining the \emph{Uniform Local Asymptotic Normality} (ULAN) property of the concentration-location FvML model (property we first establish) with Le Cam's third lemma. By doing so, we shall extend, for the FvML distribution, the ULAN property with respect to only the location parameter obtained in Ley \emph{et al.} (2013a). As we shall see, although the ULAN property does not hold for $\kappa=0$, we are nevertheless able via Le Cam's third Lemma to write down the asymptotic powers of the classical \emph{Rayleigh test} for uniformity, which is nothing but the score test for uniformity against FvML alternatives.

The rest of the paper is organized as follows. In Section~\ref{ULANsec}, we establish and prove the key ingredient of our calculations, namely the ULAN property of the concentration-location FvML model. In Sections~\ref{onesamp} and \ref{multisamp}, we write out the locally and asymptotically optimal tests for the one-sample and multi-sample problems, respectively. By construction, these coincide with the score tests proposed in Watamori and Jupp (2005) which themselves are asymptotically equivalent (the difference is $o_{\rm P}(1)$) to the likelihood ratio tests under the null (and therefore also under contiguous alternatives). In each section, we then study the asymptotic behavior of these most efficient tests under local alternatives and provide the announced expressions of their powers. In Section~\ref{onesamp}, we also study the asymptotic properties and powers (against FvML alternatives) of the famous  test for uniformity over $\mathcal{S}^{k-1}$ proposed in Rayleigh~(1919). The finite-sample powers of the tests are investigated in Section \ref{simussec} by Monte Carlo simulations, and an appendix collects the technical proofs.    

\section{The ULAN property of the concentration-location FvML model}\label{ULANsec}

Let the data points $\Xb_1,\ldots,\Xb_n$ be \mbox{i.i.d.} with common FvML density~(\ref{density}). We denote their joint distribution by ${\rm P}_{\varthetab}\n$ with $\varthetab:=(\kappa, \thetab\pr)\pr \in \R^+_0\times\mathcal {S}^{k-1}$. As announced in the Introduction, the objective of this section is to establish and prove the ULAN property of the sequence of FvML experiments $\left\{{\rm P}_{\varthetab} \n, \varthetab \in \R^+_0\times\mathcal {S}^{k-1} \right\}$. Such a sequence  is ULAN (with contiguity rate $n^{-1/2}$) if, for any sequence $\varthetab\n \in \R^+_0\times\mathcal {S}^{k-1}$ such that $\varthetab\n- \varthetab=O(n^{-1/2})$, the likelihood ratio between ${\rm P}_{\varthetab\n+n^{-1/2} \taub\n}\n$ and ${\rm P}_{\varthetab\n}\n$ allows a specific form of (probabilistic) Taylor expansion as a function of the perturbation $\taub\n \in\R\times\R^k$. In view of the curved parameter set $\R^+_0\times\mathcal {S}^{k-1}$, it is clear that the local perturbations $\taub\n$ cannot be chosen without care, as they need to satisfy that $\varthetab\n+n^{-1/2} \taub\n$ remains in $\R^+_0\times\mathcal {S}^{k-1}$. Hence, writing $\taub\n$ as $(c\n,(\tb\n)\pr)\pr$ with $c\n\in\R$ and $\tb\n\in\R^k$, we have the conditions 
\begin{equation}\label{kapcond}
\kappa\n+n^{-1/2} c\n>0
\end{equation}
and 
\begin{align} \label{loccond} 
0 &= (\thetab\n+ n^{-1/2}\tb\n)\pr ({\thetab\n}+ n^{-1/2}\tb\n) -1 \nonumber \\ 
&=  2n^{-1/2} {\thetab\n}\pr \tb\n + n^{-1} (\tb\n)\pr
\tb\n.\end{align}
The second condition thus means that the perturbation $\tb\n$ must
belong, up to a $o(n^{-1/2})$ quantity, to the tangent space to
$\mathcal{S}^{k-1}$ at $\thetab\n$. 

In order to ease readability, we introduce some notations. It can be shown that the projections $\Xb_1\pr \thetab, \ldots, \Xb_n\pr \thetab$ are \mbox{i.i.d.} with common density $\tilde{f}_\kappa(t)$ proportional to $\exp(\kappa{t}) (1-t^2)^{(k-3)/2}$ for $t\in[-1,1]$. Straightforward calculations then reveal that 
$${\rm E}[\Xb_i]={\rm E}[\Xb_i\pr\thetab] \thetab=:A_k(\kappa) \thetab =\left( \frac{\int_{-1}^1 t e^{\kappa t} (1-t^2)^{\frac{k-3}{2}}dt}{\int_{-1}^1 e^{\kappa t} (1-t^2)^{\frac{k-3}{2}}dt} \right) \thetab,$$
showing that the parameter $\kappa$ is identified via the function $A_k(\cdot)$. Note in passing that $A_k(\cdot)=I_{k/2}(\cdot)/I_{k/2-1}(\cdot)$.  Similar manipulations yield 
$${\rm Var}[\Xb_i\pr \thetab]=A_k\pr(\kappa)=1-\frac{k-1}{\kappa} A_k(\kappa)-(A_k(\kappa))^2;$$
see Watson~(1983) for more details on these results. We are now ready to state the ULAN property of the FvML concentration-location model.

\begin{Prop}\label{ULAN} The family $\left\{ {\rm P}_{{\varthetab}}\n \;  \vert \;\varthetab\in \R^+_0\times\mathcal{S}^{k-1} \right\}$ is ULAN; more precisely,  for any sequence $\varthetab\n \in \R^+_0\times\mathcal{S}^{k-1}$ such that $\varthetab\n- \varthetab=O(n^{-1/2})$ and any bounded sequence $\taub\n\in\R\times\R^k$ subjected to the conditions~(\ref{kapcond}) and~(\ref{loccond}), we have
$$\log \left( \frac{{\rm dP}\n_{{\varthetab\n} + n^{-1/2} {\taub}\n}}{{\rm dP}_{\varthetab\n}\n}\right)=  ({\taub}\n)\pr \Deltab_{\varthetab\n}\n - \frac{1}{2}({\taub}\n)\pr \Gamb_{\varthetab} {\taub}\n +o_{\rm P}(1)$$
and $\Deltab_{\varthetab\n}\n\stackrel{\mathcal{D}}{\rightarrow}\mathcal{N}_{k+1}(\zerob,\Gamb_{\varthetab})$ under ${\rm P}_{{\varthetab\n}}\n$ as $\ny$. The central sequence $\Deltab_{{\varthetab}}\n:= \left(\left(\Delta_{{\varthetab}}^{({\rm I})(n)}\right)\pr, \left(\Deltab_{\varthetab}^{({\rm II})(n)}\right)\pr \right)\pr$ is defined by
\begin{equation*}\label{centrseqbis}\Delta_{{\varthetab}}^{({\rm I})(n)}:= n^{-1/2} \sum_{i=1}^{n}( \Xb_i\pr \thetab- A_k(\kappa)) \end{equation*}
and
\begin{equation*}\label{centrseq}\Deltab_{{\varthetab}}^{({\rm II})(n)}:= \kappa n^{-1/2} \sum_{i=1}^{n}  (1-(\Xb_i \pr \thetab)^2)^{1/2} \Sb_\thetab(\Xb_i) \end{equation*}
with $\Sb_\thetab(\Xb_i):=(\Xb_i-(\Xb_i\pr\thetab)\thetab)/||\Xb_i-(\Xb_i\pr\thetab)\thetab||$.   The associated Fisher information is given by $$\Gamb_{\varthetab}:=\left( \begin{array}{cc} \Gamma_{\varthetab}^{({\rm I})} & {\bf 0} \\ {\bf 0} & \Gamb_{\varthetab}^{({\rm II})} \end{array} \right),$$
where, putting $\mathcal{J}_k(\kappa):= \int_{-1}^{1} (1-u^2) \tilde{f}_\kappa(u) du$,
\begin{equation*} \label{Fish} \Gamma_{\varthetab}^{({\rm I})}:=1-\frac{k-1}{\kappa} A_k(\kappa)-(A_k(\kappa))^2 \quad {\rm and} \quad \Gamb_{\varthetab}^{({\rm II})}:= \frac{\kappa^2 \mathcal{J}_k(\kappa)}{k-1} ({\bf I}_{k}- \thetab\thetab\pr). \end{equation*}
\end{Prop} \vspace{4mm}

This proposition constitutes, for FvML distributions on the hyperspheres $\mathcal{S}^{k-1}$, the desired extension (for FvML distributions) of Proposition~2.2 in Ley \emph{et al.}~Ê(2013a) where only the location parameter $\thetab$ was taken into account. Note the diagonal structure of the information matrix; it is the structural reason why replacing $\thetab$ by a root-$n$ consistent estimator has no asymptotic effect on inferential procedures focussing on $\kappa$. \vspace{4mm}

\noindent {\bf Proof.} We clearly need to circumvent the curved nature of the parameter space $\R^+_0\times\mathcal{S}^{k-1}$, more precisely of $\mathcal{S}^{k-1}$. Fortunately, this has been achieved in Ley \emph{et al.}~(2013a) by proving ULAN rather for the spherical coordinates $\thetab=h(\etab)$ for some locally full rank chart $h$ and then returning (via a result in Hallin \emph{et al.}~2010) to the initial $\thetab$-parameterization. Thus, thanks to the developments in Ley \emph{et al.}~(2013a), all we need to do here is to prove ULAN with respect to the ``linear'' parameters $\kappa$ and $\etab$ with $\etab\in\R^{k-1}$.

Our proof of that ULAN result relies on Lemma 1 of Swensen (1985)--more precisely, on its extension in Garel and Hallin (1995). Seven conditions need to be satisfied; we leave them to the reader, as they are easily obtained once it is proved that the mapping $(\kappa,\etab)\mapsto f_{\kappa,h(\etab)}^{1/2}(\xb)$ is differentiable in quadratic mean. The latter differentiability in quadratic mean spells out as
$$\int_{\mathcal{S}^{k-1}} \left( f_{\kappa+s,h(\etab+{\bf e})}^{1/2}(\xb)- f_{\kappa,h(\etab)}^{1/2}(\xb)- (s,{\bf e}\pr)\left(
\begin{array}{c}
\partial_\kappa f_{\kappa ,h(\etab)}^{1/2}(\xb)\\
{\rm grad}_{\etab}  f_{\kappa ,h(\etab)}^{1/2}(\xb)
\end{array}
\right)\right)^2 \; d\sigma(\xb)= o\left(\left\|\begin{array}{c}s\\{\bf e}\end{array}\right\|^2\right)$$
for $s\in\R$ and ${\bf e}\in\R^{k-1}$ such that $\kappa+s>0$ and $h(\etab+{\bf e})\in\mathcal{S}^{k-1}$. This result holds true once we have demonstrated the following three equalities:
\begin{enumerate}
\item[(i)] $\int_{\mathcal{S}^{k-1}} \left( f_{\kappa ,h(\etab+{\bf e})}^{1/2}(\xb)- f_{\kappa ,h(\etab)}^{1/2}(\xb)- ({\rm grad}_{\etab}  f_{\kappa ,h(\etab)}^{1/2}(\xb))\pr {\bf e}\right)^2 \; d\sigma(\xb)= o(\| {\bf e} \|^2)$;
\item[(ii)] $\int_{\mathcal{S}^{k-1}} \left( f_{\kappa +s,h(\etab+{\bf e})}^{1/2}(\xb)- f_{\kappa ,h(\etab+{\bf e})}^{1/2}(\xb)-s\partial_\kappa f_{\kappa ,h(\etab+{\bf e})}^{1/2}(\xb)\right)^2 \; d\sigma(\xb)= o(s^2)$, and
\item[(iii)]  $\int_{\mathcal{S}^{k-1}} \left(\partial_\kappa f_{\kappa ,h(\etab+{\bf e})}^{1/2}(\xb)-\partial_\kappa f_{\kappa ,h(\etab)}^{1/2}(\xb)\right)^2 \; d\sigma(\xb)= o(1)$.
\end{enumerate}

Point (i) has been obtained in Ley \emph{et al.} (2013a). Now for Point (ii), first note that, letting $\tilde\thetab:=h(\etab+{\bf e}) \in \mathcal{S}^{k-1}$, we have 
\begin{eqnarray*} \partial_\kappa f^{1/2}_{\kappa,{\tilde \thetab}} &=&  \frac{c_{k, \kappa}^{1/2}\, \xb\pr{\tilde \thetab}}{2} {\rm exp}(\kappa \xb\pr{\tilde \thetab}/2)+  (\partial_\kappa c_{k, \kappa}^{1/2}) {\rm exp}(\kappa \xb\pr{\tilde \thetab}/2).
\end{eqnarray*}
Therefore, the integral of Point (ii) can be rewritten as
$$
\int_{\mathcal{S}^{k-1}} {\rm exp}(\kappa \xb\pr{\tilde \thetab}) \left(c^{1/2}_{k, \kappa+s}{\rm exp}(s \xb\pr{\tilde \thetab}/2)-c^{1/2}_{k, \kappa}
-s \left(c_{k, \kappa}^{1/2} \frac{\xb\pr{\tilde \thetab}}{2} + \partial_\kappa c_{k, \kappa}^{1/2} \right)
\right)^2 \; d\sigma(\xb).
$$
This integral can be bounded by $c_1S_1+c_2S_2+c_3s^2 S_3$, where
$$S_1:= \int c_{k, \kappa+s} \; {\rm exp}(\kappa \xb\pr{\tilde \thetab}) \; \left({\rm exp}(s \xb\pr{\tilde \thetab}/2)-1- \frac{s \xb\pr \tilde{\thetab}}{2} \right)^2 d\sigma(\xb),$$
$$S_2:=\int \left(c_{k, \kappa+s}^{1/2}-c_{k, \kappa}^{1/2}-s\partial_\kappa c_{k, \kappa}^{1/2} \right)^2 \; {\rm exp}(\kappa \xb\pr{\tilde \thetab}) \; d\sigma(\xb),$$
and 
$$S_3:=\int \left(c_{k, \kappa+s}^{1/2}-c_{k, \kappa}^{1/2}\right)^2 \;  \left(\frac{\xb\pr{\tilde \thetab}}{2}\right)^2 {\rm exp}(\kappa \xb\pr{\tilde \thetab}) \; d\sigma(\xb).$$
Since both ${\rm exp}(\kappa \xb\pr{\tilde \thetab})$ and $ (\xb\pr{\tilde \thetab})^2 \; {\rm exp}(\kappa \xb\pr{\tilde \thetab})$ are obviously integrable on $\mathcal{S}^{k-1}$, 
it follows from the derivability of the mapping $\kappa \mapsto c_{k, \kappa}^{1/2}$ that $S_3$ is o(1) and that $S_2$ is $o(s^2)$. Now, the derivability of the mapping $t \mapsto {\rm exp}(t)$ at $t=0$ combined with Lebesgue's dominated convergence theorem directly entails that $S_1$ is $o(s^2)$. Putting the ends together, we have proved Point (ii).

Finally, Point (iii) follows along the same lines since all quantities involved are differentiable and square-integrable. This concludes the proof. \cqfd

\section{One-sample tests on the concentration parameter} \label{onesamp}

In this section, our focus lies on the one-sample testing problem $\mathcal{H}_0^{\kappa_0}: \kappa= \kappa_0$ for some fixed $\kappa_0 >0$ versus $\mathcal{H}_1^{\kappa_0}: \kappa\neq \kappa_0$ (Section~\ref{ones}) and on the Rayleigh~(1919) tests of uniformity (Section~\ref{Ray}). In each case, we analyze the most efficient tests (Watamori and Jupp~2005 and Rayleigh~1919, respectively) in the light of the Le Cam framework. Besides interesting optimality issues, this will allow us to obtain expressions for the powers of the considered tests under local alternatives.

\subsection{The one-sample score tests of Watamori and Jupp~(2005)}\label{ones}

Based on the ULAN property in Proposition~\ref{ULAN}, the Le Cam asymptotic theory (see Le Cam~1986) paves the way towards constructing locally and asymptotically optimal tests. The optimality appearing in this section is the so-called \emph{maximin} optimality. A test $\phi^*$ is called maximin in the class $\mathcal{C}_\alpha$ of level-$\alpha$ tests  for $\mathcal{H}_0$ against $\mathcal{H}_1$ if (i) $\phi^*$ has level $\alpha$ and (ii) the power of $\phi^*$ is such that
$$\inf_{{\rm P} \in \mathcal{H}_1}{\rm E}_{{\rm P}}[\phi^*] \geq \sup_{\phi \in \mathcal{C}_\alpha} \inf_{{\rm P} \in \mathcal{H}_1} {\rm E}_{\rm P}[\phi].$$
Since $\kappa$ is the parameter of interest, locally and asymptotically optimal tests for $\mathcal{H}_0^{\kappa_0}$ are built upon $\Delta_{\varthetab}^{({\rm I})(n)}$, the $\kappa$-part of the central sequence; see Le Cam~(1986), Section 11.9, for details. More concretely, a locally and asymptotically maximin test rejects $\mathcal{H}_{0}^{\kappa_0}$ at asymptotic level $\alpha$ whenever
\begin{eqnarray*}
Q_{\kappa_0}\n (\thetab)&:=&\frac{(\Delta_{\varthetab}^{({\rm I})(n)})^{2}}{\Gamma_{\varthetab}^{({\rm I})}} \nonumber \\
&=& \frac{\left(\sum_{i=1}^{n}(\Xb_i\pr \thetab- A_k^{-1}(\kappa_0))\right)^2}{n(1-\frac{k-1}{\kappa_0} A_k(\kappa_0)-(A_k(\kappa_0))^2)}
\end{eqnarray*}
exceeds the $\alpha$-upper quantile of the chi-square distribution with 1 degree of freedom. Unfortunately,  $Q_{\kappa_0} \n(\thetab)$ is not (yet) a genuine test statistic since it still depends on the unknown value of $\thetab$. This problem can be solved by replacing $\thetab$ with a root-$n$ consistent estimator $\hat\thetab\n$ in the central sequence $\Delta_{\varthetab}^{({\rm I})(n)}$, whilst, of course, paying attention to the asymptotic effects of such a substitution. It is here that the ULAN property of the concentration-location FvML model comes in handy. Indeed, it directly entails (see again Le Cam~1986) that the FvML model is locally and asymptotically linear in the sense that
\begin{equation} \label{alinearity}
\Deltab_{\varthetab+n^{-1/2} \taub\n}\n-\Deltab_{\varthetab}\n=\Gamb_{\varthetab} \taub\n + o_{\rm P}(1)
\end{equation}
under ${\rm P}_{{\varthetab}}\n$ as $\ny$. Of course, the aim consists in using $\taub\n=(c\n,(\tb\n)\pr)\pr$ with $\tb\n=n^{1/2}(\hat\thetab\n-\thetab)$ which satisfies condition~(\ref{loccond}); controlling this replacement however is not straightforward and requires a formal proof. The matters are simplified by the (already discussed) block-diagonality of the Fisher information matrix, which implies that the $\kappa$-part of the central sequence is not influenced by a local perturbation of $\thetab$ (similarly, the $\thetab$-part of the central sequence is not influenced by a local perturbation of $\kappa$, hence the results in Ley \emph{et al.}~2013a for the FvML case can be extended by estimating $\kappa$). Hence, since our focus lies on $\Delta_{\varthetab}^{({\rm I})(n)}$ and since, under the null hypothesis, $\kappa$ is fixed to $\kappa_0$, we only need to show by having recourse to the asymptotic linearity property~(\ref{alinearity}) for $\Delta_{\varthetab}^{({\rm I})(n)}$ that a replacement of $\thetab$ with a root-$n$ consistent estimator $\hat\thetab\n$ (e.g., the sample spherical mean $\bar\Xb/||\bar\Xb||$ with $\bar\Xb=n^{-1}\sum_{i=1}^n\Xb_i$) has no asymptotic impact on $\Delta_{\varthetab}^{({\rm I})(n)}$, which is achieved in the following proposition.

\begin{Prop}\label{alinea} Let $\hat{\thetab}\n \in {\mathcal S}^{k-1}$ be a root-$n$ consistent estimator of $\thetab$ under ${\rm P}_\varthetab\n$. Then, letting $T\n_\thetab:=n^{-1/2} \sum_{i=1}^{n} \Xb_i\pr \thetab$, we have that $T\n_{\hat \thetab\n}-T\n_{\thetab}$ is $o_{\rm P}(1)$ under ${\rm P}_\varthetab\n$ as $\ny$. 
\end{Prop}

See the appendix for the proof. The resulting locally and asymptotically maximin test for testing $\mathcal{H}_{0}^{\kappa_0}$ rejects the null (at asymptotic level $\alpha$) when
\begin{eqnarray*} \label{statkappanot} Q_{\kappa_0}\n:= \frac{\left(\sum_{i=1}^{n} \Xb_i\pr \hat\thetab- A_k^{-1}(\kappa_0)\right)^2}{n(1-\frac{k-1}{\kappa_0} A_k(\kappa_0)-(A_k(\kappa_0))^2)}
\end{eqnarray*}
exceeds the $\alpha$-upper quantile of the chi-square distribution with 1 degree of freedom. The test statistic $Q_{\kappa_0}\n$ coincides with the score test proposed in Watamori and Jupp (2005). This, in passing, shows the local and asymptotic optimality property of the latter. The following result characterizes the asymptotic properties of $Q_{\kappa_0}\n$.

\begin{Prop}\label{chideuxone}
We have that
\begin{itemize}
\item[(i)] $Q_{\kappa_0}\n$ is asymptotically chi-square with $1$ degree of freedom under  $\cup_{\thetab\in\mathcal{S}^{k-1}}{\rm P}_{(\kappa_0,\thetab)}\n$;
\item[(ii)] $Q_{\kappa_0}\n$ is asymptotically non-central chi-square with $1$ degree of freedom and non-centrality parameter 
$ (1-\frac{k-1}{\kappa_0} A_k(\kappa_0)-(A_k(\kappa_0))^2) c^2$ under $\cup_{\thetab\in\mathcal{S}^{k-1}}{\rm P}\n_{(\kappa_0+n^{-1/2} c\n,\thetab)}$ ($c:=\lim_{\ny} c\n$ for $c\n$ satisfying condition~(\ref{kapcond}));
 \item[(iii)] the test $\phi_{\kappa_0}\n$ which rejects the null hypothesis as soon as $Q_{\kappa_0}\n$ exceeds the $\alpha$-upper quantile of the chi-square distribution with $1$ degree of freedom has 
asymptotic level $\alpha$ under $\cup_{\thetab\in\mathcal{S}^{k-1}}{\rm P}\n_{(\kappa_0,\thetab)}$ and is locally and asymptotically maximin against local alternatives of the form $\cup_{\thetab\in\mathcal{S}^{k-1}}{\rm P}\n_{(\kappa_0+n^{-1/2} c\n,\thetab)}$.
\vspace{1mm}
\end{itemize}
\end{Prop}

Proposition~\ref{chideuxone} readily follows from Proposition~\ref{ULAN}, Proposition~\ref{alinea} and the celebrated third Lemma of Le Cam, and is hence left to the reader (if unclear, see the next section where we develop this argument for the Rayleigh~1919 test of uniformity). Note that Proposition \ref{chideuxone} readily yields the announced expression for the power of $Q_{\kappa_0}\n$ under local alternatives of the form $\cup_{\thetab\in\mathcal{S}^{k-1}}{\rm P}\n_{(\kappa_0+n^{-1/2} c\n,\thetab)}$ ($c:= \lim_{\ny} c\n$): 
$$1-F_{\chi^2_1((1-\frac{k-1}{\kappa_0} A_k(\kappa_0)-(A_k(\kappa_0))^2) c^2)}(\chi^2_{1;1-\alpha}),$$
where $F_{\chi^2_\nu(z)}$ stands for the distribution function of the non-central chi-square distribution with $\nu$ degrees of freedom and with non-centrality parameter $z$ and $\chi^2_{\nu;1-\alpha}$ represents the $\alpha$-upper quantile of the (central) chi-square distribution with $\nu$ degrees of freedom. 

\subsection{The Rayleigh~(1919) tests of uniformity}\label{Ray}

Let us now come to the Rayleigh test of uniformity. Within the FvML family, the boundary distribution obtained when $\kappa = 0$ is the uniform distribution. Unfortunately, the ULAN property of Proposition~\ref{ULAN} does not hold for $\kappa=0$ (\emph{inter alia} because the location $\thetab$ is not identified under the null of uniformity). Nevertheless, we show in this subsection that a study of the asymptotic local powers of the  classical Rayleigh~(1919) test which rejects the null hypothesis of uniformity $\mathcal{H}_0^{\rm unif}$ at asymptotic nominal level $\alpha$ when
 \begin{eqnarray}\label{Rayl} {\rm Q}_{\rm unif}\n := kn \| \bar{\Xb}\|^2 > \chi^2_{k; 1- \alpha}
 \end{eqnarray}
 can be performed using the Third Le Cam Lemma. It follows from (\ref{Rayl}) that in order to obtain local powers of $ {\rm Q}_{\rm unif}\n$ we have to study the asymptotic behavior of ${\bf T}\n:=n^{1/2} \bar{\Xb}$ under local FvML alternatives since ${\rm Q}_{\rm unif}\n=k ({\bf T}\n)\pr{\bf T}\n $.
 First, let $$\Lambda\n:= \log \left(\frac{{\rm dP}\n_{(n^{-1/2}c\n,\thetab)}}{{\rm dP}\n_{\rm unif}} \right)$$ stand for the log-likelihood ratio between a FvML distribution with parameters $(n^{-1/2}c\n, \thetab)$ and the uniform distribution on $\mathcal{S}^{k-1}$. Both distributions are clearly contiguous. Very simple computations yield
$$\Lambda\n= n^{-1/2} c\n \sum_{i=1}^n \Xb_i\pr \thetab + C_{k,c\n}\n$$
for some constant $C_{k,c\n}\n$ which is $o(1)$ as $\ny$ under ${{\rm P}\n_{\rm unif}}$. The multivariate central limit theorem  directly entails that the limiting distribution of $\left(({\bf T}\n)\pr, \Lambda\n\right)\pr$ is a $(k+1)$-variate Gaussian distribution with mean zero and covariance matrix ($c:= \lim_{\ny} c\n$)
$$\left(\begin{array}{cc} k^{-1} {\bf I}_k & ck^{-1} \thetab \\ c k^{-1} \thetab\pr & c^2k^{-1} \end{array}\right)$$
under ${\rm P}\n_{\rm unif}$ as $\ny$ (this holds for any fixed $\thetab \in \mathcal{S}^{k-1}$). Then the third Le Cam Lemma entails that the limiting distribution of ${\bf T}\n$ is a $k$-variate Gaussian distribution with mean $ck^{-1} \thetab$ and covariance matrix $ k^{-1} {\bf I}_k$ under ${\rm P}\n_{(n^{-1/2}c\n, \thetab)}$ as $\ny$. Wrapping up, we obtain the following result.
 
 \begin{Prop}\label{chideuxthree}
We have that
\begin{itemize}
\item[(i)] (Rayleigh 1919) $ {\rm Q}_{\rm unif}\n$ is asymptotically chi-square with $k$ degrees of freedom under $\mathcal{H}_{0}^{{\rm unif}}$;
\item[(ii)] $ {\rm Q}_{\rm unif}\n$ is asymptotically non-central chi-square with $k$ degrees of freedom and non-centrality parameter 
$ c^2/k$ under $\cup_{\thetab\in\mathcal{S}^{k-1}}{\rm P}\n_{(n^{-1/2} c\n,\thetab)}$ ($c:= \lim_{\ny} c\n$ for $c\n$ satisfying condition~(\ref{kapcond}));
\end{itemize}  
\end{Prop}

\begin{figure}[htbp!] 
\begin{center}
\includegraphics[scale=0.5]{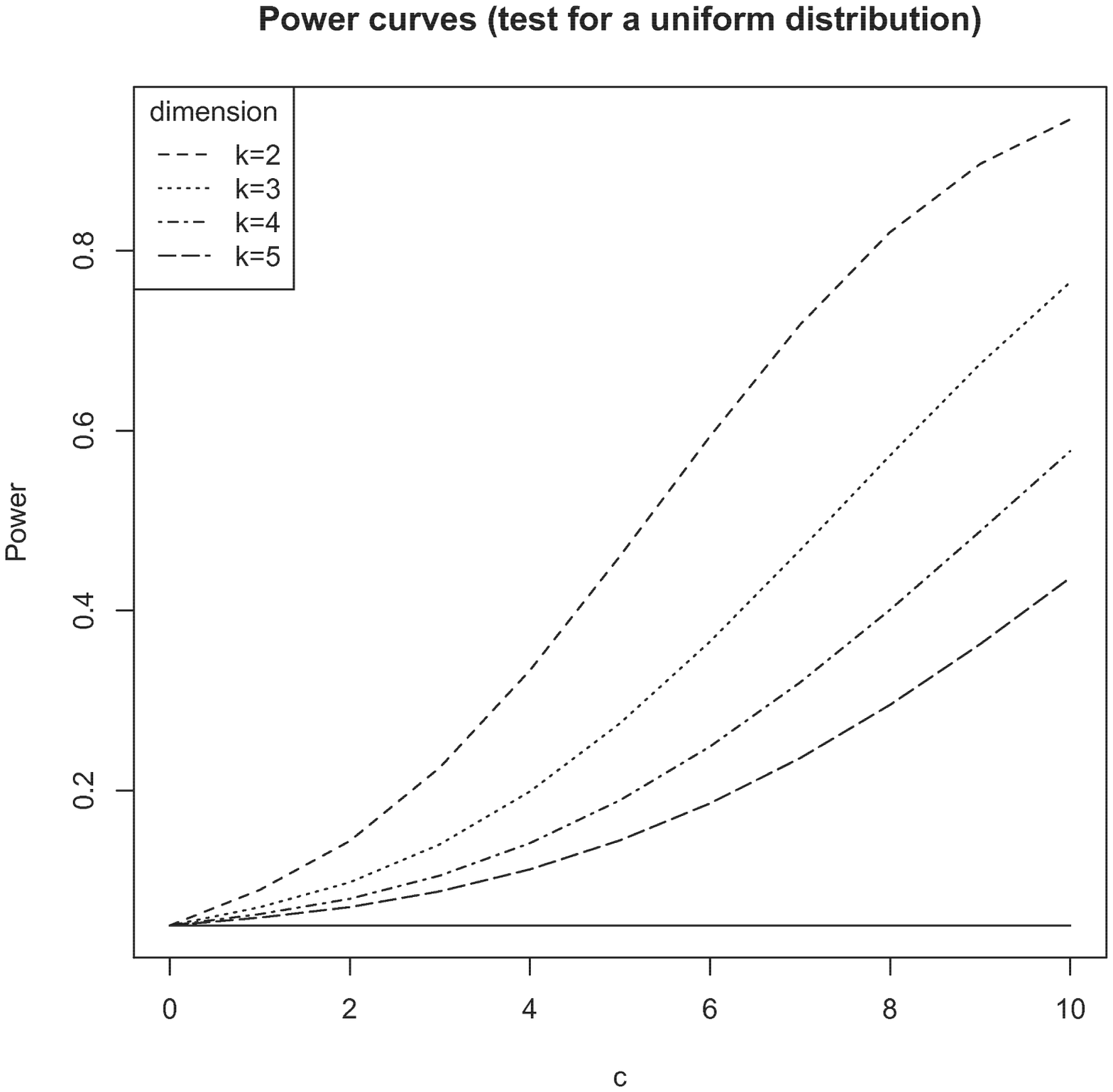}
\caption{\small Power curve (the local alternative is ${\rm P}\n_{(n^{-1/2} (c/2),\thetab)}$) of $\phi\n_{\rm unif}$ for dimensions $k=2, \ldots,5$.}
\label{curve}
\end{center}
\end{figure}

 The power of the Rayleigh test under local alternatives of the form $\cup_{\thetab\in\mathcal{S}^{k-1}}{\rm P}\n_{(n^{-1/2} c\n ,\thetab)}$ is given by ($c:= \lim_{\ny} c\n$)
 $$1-F_{\chi^2_k(c^2/k)}(\chi^2_{k;1-\alpha}).$$  Figure \ref{curve} right below shows power curves of the Rayleigh test for different values of the dimension $k$ against local FvML alternatives. Note that the power of the Rayleigh test decreases as the dimension $k$ increases.

\section{Multi-sample tests on the equality of concentrations} \label{multisamp}

In this section, our focus lies on the multi-sample testing problem $\mathcal{H}_0^{\rm Hom}: \kappa_1= \ldots= \kappa_m$ for $m\geq2$ and $\kappa_1,\ldots,\kappa_m>0$ versus $\mathcal{H}_1^{\rm Hom}: \exists 1\leq i\neq j\leq m\,\, \kappa_i\neq \kappa_j$. In other words, we are dealing with 
$m(\geq2)$ samples of \mbox{i.i.d.} data points $\Xb_{i1}, \ldots, \Xb_{in_i}$ with common FvML distribution with concentration $\kappa_i$  and location $\thetab_i$ for all $i=1,\ldots,m$, and we are interested in determining whether or not these $m$ samples have the same concentration parameters, \emph{without} assuming equality of the mean directions $\thetab_i$. As in the previous section, our way of proceeding consists in ``re-discovering'' the score tests of Watamori and Jupp~(2005) thanks to the ULAN property (which we shall adapt to the multi-sample case) and then unveiling new asymptotic results for these tests. 

Let us assume that the  samples $(\Xb_{i1}, \ldots, \Xb_{in_{i}})$, $i=1, \ldots, m$, are  mutually independent samples of \mbox{i.i.d.} random vectors; as already mentioned above, the $n_i$ observations $\Xb_{ij}$, $j=1, \ldots, n_i,$ in sample $i$ are \mbox{i.i.d.} with common FvML density with concentration $\kappa_i$ and location $\thetab_i$. We denote this time by ${\rm P}_{\varthetab^{(m)}}\n$ the joint distribution of $(\Xb_{11}, \ldots, \Xb_{mn_{m}})$, with $\varthetab^{(m)}:=(\kappa_1, \ldots, \kappa_m, \thetab_1\pr, \ldots, \thetab_m\pr)\pr \in  (\R^{+}_0)^m\times (\mathcal{S}^{k-1})^m $. In order to be able to state our results, we need to impose a
certain amount of control on the respective sample sizes $n_i$, $i=1, \ldots,
m$. This is   achieved via the following 

\

\noindent {\sc Assumption A}. Let $n=\sum_{i=1}^m n_i$. For all $i=1, \ldots, m$, the ratio
$r_i\n:=n_i/n$ converges to a non-zero constant $r_i$ as $\ny$. 

\

\noindent A direct consequence of Assumption~A is that the specific
sizes $n_i$ become somehow irrelevant; hence,  in what
follows, we simply use the superscript $^{(n)}$  for the different
quantities at play and do not specify whether they are associated with
a given $n_i$.  Now, let ${\rm diag}({\bf A}_1, \ldots, {\bf A}_m)$ stand
for the block-diagonal matrix with blocks ${\bf A}_1,
\ldots, {\bf A}_m$,  and use the notation 
$  \nub\n:={\rm diag}(\nub_1\n, \nub_2\n)$, where $\nub_1\n:={\rm diag}((r_1\n)^{-1/2}, \ldots, (r_m\n)^{-1/2})$ and $\nub_2\n:={\rm   diag}((r_1\n)^{-1/2} {\bf I}_k, \ldots, (r_m\n)^{-1/2} {\bf I}_k)
$. As in the one-sample case, we only consider perturbations $\taub\n=(c_1\n,\ldots, c_m\n,({\tb}_1\n)\pr, \ldots, ({\tb}_m\n)\pr)\pr\in\R^m\times(\R^k)^m$ such that, for any $\varthetab^{(m)} \in (\R^{+}_0)^m\times (\mathcal{S}^{k-1})^m$, $\varthetab^{(m)}+n^{-1/2}\nub\n \taub\n$ remains in $(\R^{+}_0)^m\times (\mathcal{S}^{k-1})^m$ (this is simply an adaptation of the conditions~(\ref{kapcond}) and~(\ref{loccond})). This readily leads us to the following multi-sample version of Proposition~\ref{ULAN}, whose straightforward proof is omitted.

\begin{Prop}\label{ULANbis} Let Assumption A hold. Then the family $\left\{ {\rm P}_{{\varthetab}^{(m)}}\n \;  \vert \;\varthetab^{(m)}\in (\R^{+}_0)^m\times (\mathcal{S}^{k-1})^m \right\}$ is ULAN; more precisely, for any sequence $(\varthetab^{(m)})\n\in(\R^{+}_0)^m\times (\mathcal{S}^{k-1})^m$ such that $(\varthetab^{(m)})\n-\varthetab^{(m)}=O(n^{-1/2})$ and any bounded sequence $\taub\n$ as described just before,
$$\log \left( \frac{{\rm dP}\n_{(\varthetab^{(m)})\n + n^{-1/2} \nub\n {\taub}\n}}{{\rm dP}\n_{(\varthetab^{(m)})\n}}\right)=  ({\taub}\n)\pr \Deltab_{(\varthetab^{(m)})\n}\n - \frac{1}{2}({\taub}\n)\pr \Gamb_{\varthetab^{(m)}} {\taub}\n +o_{\rm P}(1)$$
and $\Deltab_{(\varthetab^{(m)})\n}\n\stackrel{\mathcal{D}}{\rightarrow}\mathcal{N}_{m(k+1)}(\zerob,\Gamb_{\varthetab^{(m)}})$ under ${\rm P}\n_{{\varthetab}^{(m)}}$ as $\ny$. The central sequence \\$\Deltab_{{\varthetab}^{(m)}}\n \! \! := \! \! \left(\left(\Deltab_{{\varthetab}^{(m)}}^{({\rm I})(n)}\right)\pr, \left(\Deltab_{{\varthetab}^{(m)}}^{({\rm II})(n)}\right)\pr\right)\pr$, where $\Deltab_{{\varthetab}^{(m)}}^{({\rm I})(n)} \!:= \!  \!  \left(\left(\Delta_{{\varthetab}^{(m)}}^{({\rm I},1)(n)}\right)\pr, \ldots, \left(\Delta_{{\varthetab}^{(m)}}^{({\rm I},m) (n)}\right)\pr \right)\pr$ and $\Deltab_{{\varthetab}^{(m)}}^{({\rm II})(n)}:= \left(\left(\Deltab_{{\varthetab}^{(m)}}^{({\rm II},1)(n)}\right)\pr, \ldots, \left(\Deltab_{{\varthetab}^{(m)}}^{({\rm II},m)(n)} \right)\pr\right)\pr$, is defined by
 \begin{equation*}\label{centrseqbis}\Delta_{{\varthetab}^{(m)}}^{({\rm I},i)(n)}:= n_i^{-1/2} \sum_{j=1}^{n_i}(\Xb_{ij}\pr \thetab_i- A_k(\kappa_i))\end{equation*}
and
\begin{equation*}\label{centrseq}\Deltab_{{\varthetab}^{(m)}}^{({\rm II},i)(n)}:= \kappa_i n_i^{-1/2} \sum_{j=1}^{n_i}  (1-(\Xb_{ij} \pr \thetab_i)^2)^{1/2} \Sb_{\thetab_i}(\Xb_{ij})\end{equation*}
for all $i=1,\ldots,m$, with $\Sb_\thetab(\Xb_{ij}):=(\Xb_{ij}-(\Xb_{ij}\pr\thetab_i)\thetab_i)/||\Xb_{ij}-(\Xb_{ij}\pr\thetab_i)\thetab_i||$.   The associated Fisher information is given by $\Gamb_{\varthetab^{(m)}}:={\rm diag}(\Gamb_{\varthetab^{(m)}}^{({\rm I})}, \Gamb_{\varthetab^{(m)}}^{({\rm II})})$, where  $\Gamb_{\varthetab^{(m)}}^{({\rm I})}:={\rm diag}(\Gamma_{\varthetab^{(m)}}^{({\rm I},1)}, \ldots, \Gamma_{\varthetab^{(m)}}^{({\rm I}, m)})$ with $ \Gamma_{\varthetab^{(m)}}^{({\rm I},i)}:=1-\frac{k-1}{\kappa_i} A_k(\kappa_i)-(A_k(\kappa_i))^2$ for all $i=1,\ldots,m$ and where
 $\Gamb_{\varthetab^{(m)}}^{({\rm II})}:= {\rm diag}(\Gamb_{\varthetab^{(m)}}^{({\rm II},1)}, \ldots, \Gamb_{\varthetab^{(m)}}^{({\rm II},m)})$ with
 \begin{equation*} \label{Fish1} \Gamb_{\varthetab^{(m)}}^{({\rm II},i)}:= \frac{\kappa_i^2 \mathcal{J}_k(\kappa_i)}{k-1} ({\bf I}_{k}- \thetab_i\thetab_i\pr)
 \end{equation*}
for all $i=1, \ldots, m.$
\end{Prop}

As for the one-sample case in Section~\ref{onesamp}, we use the ULAN property to construct a locally and asymptotically optimal test for the homogeneity of concentrations. Here, the underpinning optimality concept provides the so-called  \emph{most stringent} test for $\mathcal{H}_0^{\rm Hom}$. A  test $\phi^*$ is called most stringent in the class of level-$\alpha$ tests $\mathcal{C}_\alpha$ for testing $\mathcal{H}_0$ against
$\mathcal{H}_1$ if (i) $\phi^*$ has level $\alpha$ and (ii) is such that 
$$
  \sup_{{\rm P}\in\mathcal{H}_1}r_{\phi^*}({\rm P})\leq\sup_{{\rm P}\in\mathcal{H}_1}r_\phi({\rm P})\quad\forall\phi\in\mathcal{C}_\alpha,
$$
where $r_{\phi_0}({\rm P})$ stands for the \emph{regret} of the test
$\phi_0$ under ${\rm P}\in\mathcal{H}_1$ defined as $r_{\phi_0}({\rm
  P}):=\left[\sup_{\phi\in\mathcal{C}_\alpha}{\rm E}_{\rm
    P}[\phi]\right]-{\rm E}_{\rm P}[\phi_0]$, the deficiency in power
of $\phi_0$ under ${\rm P}$ compared to the highest possible (for
tests belonging to $\mathcal{C}_\alpha$) power under ${\rm P}$.   

Letting ${\bf 1}_m:=(1, \ldots, 1)\pr \in \R^m$, the null hypothesis $\mathcal{H}_0^{\rm Hom}$ can be rewritten as $\mathcal{H}_0^{\rm Hom}:(\kappa_1, \ldots, \kappa_m)\pr \in \mathcal{M}({\bf 1}_m)$, where $\mathcal{M}({\bf A})$ stands for the linear subspace spanned by the columns of ${\bf A}$.  Following Le Cam (1986), a locally and asymptotically most stringent test rejects the null hypothesis $\mathcal{H}_0^{\rm Hom}$ at asymptotic level $\alpha$ when (writing $\kappa$ for the common value of $\kappa_1, \ldots, \kappa_m$ under the null, $D_k:=1-\frac{k-1}{\kappa} A_k(\kappa)-(A_k(\kappa))^2$, $\Umb\n:=(\nub_1\n)^{-1} {\bf 1}_m$ and $\bar{\Xb}_i:= n_i^{-1} \sum_{j=1}^{n_i} \Xb_{ij}$)
\begin{eqnarray}
Q_{\rm Hom}\n(\varthetab^{(m)})&:=&\left(\Deltab_{{\varthetab}^{(m)}}^{({\rm I})(n)}\right)\pr \left((\Gamb_{\varthetab^{(m)}}^{({\rm I})})^{-1}- \Umb\n \left((\Umb\n)\pr \Gamb_{\varthetab^{(m)}}^{({\rm I})}\Umb\n \right)^{-1} (\Umb\n)\pr \right) \Deltab_{{\varthetab}^{(m)}}^{({\rm I})(n)} \nonumber \\
&=& D_k^{-1} \left(\sum_{i=1}^{m}n_i (\thetab_i\pr \bar{\Xb}_i)^2- \frac{1}{n} \left( \sum_{i=1}^m n_i \thetab_i\pr \bar{\Xb}_i \right)^2 \right) \nonumber
\end{eqnarray}
exceeds the $\alpha$-upper quantile of the chi-square distribution with $m-1$ degrees of freedom. As for the one-sample case, the statistic $Q_{\rm Hom}\n(\varthetab^{(m)})$ is not (yet) a genuine test statistic since it still depends on the unknown location parameters $\thetab_1, \ldots, \thetab_m$ and moreover on the quantity $D_k$. The replacement of the location parameters with root-$n$ consistent estimators (e.g., $\hat\thetab_{1}\n={\bar\Xb}_1/ \| {\bar\Xb}_1 \|, \ldots, \hat\thetab_{m}\n= {\bar\Xb}_m/ \| {\bar\Xb}_m \|$ with ${\bar\Xb}_j/ \| {\bar\Xb}_j\|$ the $j$th intra-sample spherical mean) will not have any asymptotic impact on $Q_{\rm Hom}\n(\varthetab^{(m)})$, see Proposition~\ref{alinea}. As concerns the quantity $D_k$, it can be estimated consistently by $\hat{D}_k:= 1-\frac{k-1}{\hat{\kappa}} A_k(\hat\kappa)-(A_k(\hat\kappa))^2$, where, putting $\hat{\kappa}_1=A_k^{-1}(\| \bar{\Xb}_1\|), \ldots, \hat{\kappa}_m=A_k^{-1}(\| \bar{\Xb}_m\|)$, $\hat{\kappa}:= \sum_{i=1}^m r_i\n {\hat \kappa}_i$. The resulting locally and asymptotically most stringent test  $\phi_{\rm Hom}\n$
 rejects the null hypothesis $\mathcal{H}_0^{\rm Hom}$ at asymptotic level $\alpha$ whenever 
 \begin{eqnarray*}
Q_{\rm Hom}\n&:=& \hat{D}_k^{-1} \left(\sum_{i=1}^{m}  n_i (\hat{\thetab}_i\pr \bar{\Xb}_i)^2- \frac{1}{n} \left( \sum_{i=1}^m n_i \hat{\thetab}_i\pr \bar{\Xb}_i \right)^2 \right)
\end{eqnarray*}
exceeds the $\alpha$-upper quantile of the chi-square distribution with $m-1$ degrees of freedom. 
Again, the test statistic $Q_{\rm Hom}\n$ coincides with the score test proposed in Watamori and Jupp (2005) which is therefore locally and asymptotically most stringent. The following result characterizes the asymptotic properties of  $Q_{\rm Hom}\n$ under the null and under a sequence of local alternatives.

\begin{Prop}\label{chideuxtwo}
Let Assumption A hold. We have that
\begin{itemize}
\item[(i)] $Q_{\rm Hom}\n$ is asymptotically chi-square with $m-1$ degrees of freedom under  $\mathcal{H}_0^{\rm Hom}$;
\item[(ii)] letting ${\bf c}=(c_1, \ldots, c_m):= \lim_{\ny} (c_1\n, \ldots, c_m\n)\pr$, $Q_{\rm Hom}\n$ is asymptotically non-central chi-square with $m-1$ degrees of freedom and non-centrality parameter 
\begin{equation}\label{noncentr}
D_k \left[\sum_{i=1}^m c_i^2-\left(\sum_{i=1}^m \sqrt{r_i} c_i\right)^2 \right]
\end{equation} 
under $\cup_{(\thetab_1\pr,\ldots,\thetab_m\pr)\pr\in(\mathcal{S}^{k-1})^m}\cup_{\kappa\in\R_0^+}{\rm P}\n_{(\kappa,\ldots,\kappa,\thetab_1\pr,\ldots,\thetab_m\pr)\pr+ n^{-1/2} \nub\n \taub\n}$;
 \item[(iii)] the test $\phi_{\rm Hom}\n$ which rejects the null hypothesis as soon as $Q_{\rm Hom}\n$ exceeds the $\alpha$-upper quantile of the chi-square distribution with $m-1$ degrees of freedom has 
asymptotic level $\alpha$ under $\mathcal{H}_0^{\rm Hom}$ and is locally and asymptotically most stringent against local alternatives of the form $\cup_{(\thetab_1\pr,\ldots,\thetab_m\pr)\pr\in(\mathcal{S}^{k-1})^m}\cup_{\kappa\in\R_0^+}{\rm P}\n_{(\kappa,\ldots,\kappa,\thetab_1\pr,\ldots,\thetab_m\pr)\pr+ n^{-1/2} \nub\n \taub\n}$.
\vspace{1mm}
\end{itemize}
\end{Prop}

See the appendix for a proof. Note that, when all quantities $c_i\n(r_i\n)^{-1/2}$ (and hence the limit $c_ir_i^{-1/2}$) are equal, we are still under the null; this is well translated by the fact that then the non-centrality parameter in~(\ref{noncentr}) equals zero. Proposition \ref{chideuxtwo} also readily yields the announced expression for the power of $Q_{\rm Hom}\n$ under local alternatives of the form $\cup_{(\thetab_1\pr,\ldots,\thetab_m\pr)\pr\in(\mathcal{S}^{k-1})^m}\cup_{\kappa\in\R_0^+}{\rm P}\n_{(\kappa,\ldots,\kappa,\thetab_1\pr,\ldots,\thetab_m\pr)\pr+ n^{-1/2} \nub\n \taub\n}$:
$$1-F_{\chi^2_{m-1}\left( D_k \left[\sum_{i=1}^m c_i^2-\left(\sum_{i=1}^m \sqrt{r_i} c_i\right)^2 \right] \right)}(\chi^2_{m-1;1-\alpha}).$$
 We conclude this section by attracting the reader's attention to the fact that this multi-sample problem here complements, for the FvML case, the ANOVA study in Ley \emph{et al.}~(2013b).

In the next section, we study the finite-sample powers of the tests constructed here via Monte Carlo simulations.

\section{Monte Carlo simulations}\label{simussec}

Since Watamori and Jupp~(2005) do not examine the finite-sample performances of their score tests, we will do so via a Monte Carlo study in this section. More precisely, we shall concentrate on the multi-sample case and hence complement the theoretical powers provided at the end of the previous section by a simulation study. However, before starting this analysis, we will first verify numerically the asymptotic powers obtained for the Rayleigh (1919) test.   

\subsection{Power curve of the Rayleigh (1919) test}

\begin{figure} [htbp!] 
\begin{center}
\includegraphics[scale=0.5]{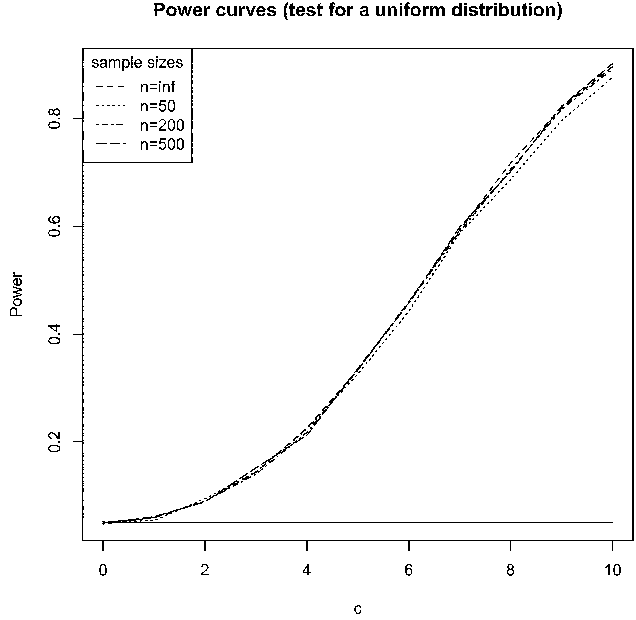}
\caption{\small Power curve of $\phi\n_{\rm unif}$ for $k=2$ and various sample sizes $n=50$, $n=100$ and $n=200$. The ``n=inf" curve is the true (theoretical) power curve (obtained using Proposition \ref{chideuxthree}).}
\label{emppoRay}
\end{center}
\end{figure}

The aim of this subsection is to corroborate Proposition \ref{chideuxthree} and the ensuing power curves by showing that empirical power curves do converge to the theoretical ones. To do so, we generated $N=5,000$ independent replications of circular FvML (hence, in fact, von Mises) random vectors 
$$
{\Xb}_{c; i}, \ \ \ c=0, \ldots, 10, \quad  {i}=1,\ldots, n,
$$
with concentration $n^{-1/2} c/2$ and location $\thetab=(1,0)\pr$. The vectors ${\Xb}_{0; i}$ represent  the null hypothesis while the vectors ${\Xb}_{c; i}$ for $c=1, \ldots, 10$ are (increasingly) under the alternative. The results using sample sizes $n=50$, $n=200$ and $n=500$ are plotted in Figure \ref{emppoRay}. They clearly confirm the theoretical power curves and hence Proposition \ref{chideuxthree}.

\subsection{Finite sample behavior of the most stringent test for the homogeneity of concentrations}

\begin{figure}[htbp!] 
\begin{center}
\includegraphics[scale=0.5]{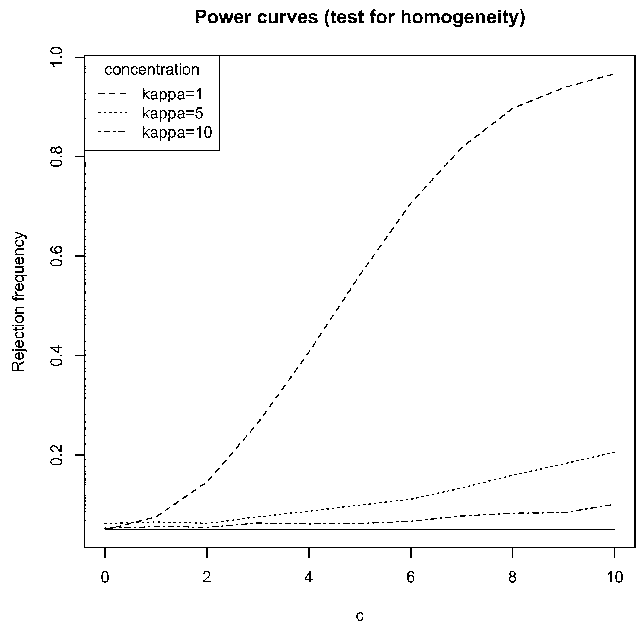}
\caption{\small Empirical power curves of $\phi\n_{\rm Hom}$ for $k=2$, various concentration parameters ($\kappa=1, 5, 10$) and sample sizes $n_1=100$ and $n_2=150$}
\label{fighom1}
\end{center}
\end{figure}

\begin{figure}[htbp!] 
\begin{center}
\includegraphics[scale=0.5]{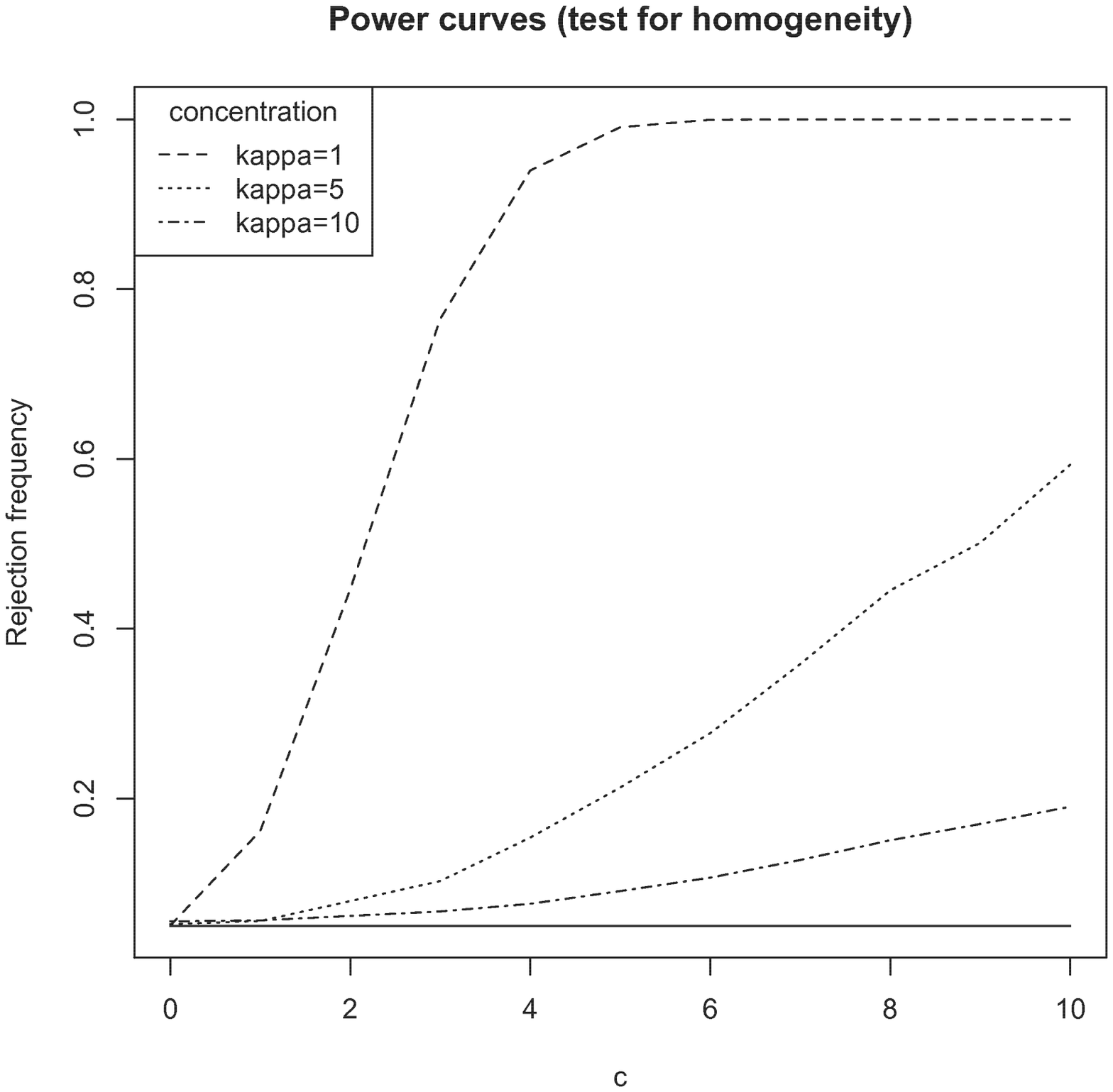}
\caption{\small Empirical power curves of $\phi\n_{\rm Hom}$ for $k=2$, various concentration parameters  ($\kappa=1, 5, 10$) and sample sizes $n_1=500$ and $n_2=500$}
\label{fighom2}
\end{center}
\end{figure}

In this subsection, we investigate the finite-sample behavior of the test $\phi_{\rm Hom}\n$ for the homogeneity of the concentrations.
We generated $N=5,000$ independent replications of three pairs ($m=2$) of mutually independent samples (we considered two designs; first with respective, and relatively small, sizes $n_1=100$ and $n_2=150$ and then with respective moderate sample sizes $n_1=500$ and $n_2=500$) of circular  random vectors
$$
{\Xb}_{\ell;1j_{1}}
 \quad {\rm and} \quad {\Xb}_{\ell c ;2j_{2}}, \quad
\ell=1, 2, 3, \ \ \ c=0, \ldots, 10, \quad  j_{i}=1,\ldots, n_{i}, \quad i=1,2,
$$
with FvML densities such that 
\begin{itemize}
\item[(i)] ${\Xb}_{1;1j_1}$ and ${\Xb}_{10;2j_2}$ have a common concentration $\kappa=1$ and locations $\thetab_1=(1,0)\pr$ and $\thetab_2=(-1,0)\pr$. Then for $c=1, \ldots, 10$,  the ${\Xb}_{1c;2j_2}$'s have concentration $1+c/10$ and still locations $\thetab_1=(1,0)\pr$ and $\thetab_2=(-1,0)\pr$.
 \item[(ii)] ${\Xb}_{2;1j_1}$ and ${\Xb}_{20;2j_2}$ have a common concentration $\kappa=5$ and locations $\thetab_1=(1,0)\pr$ and $\thetab_2=(-1,0)\pr$. Then for $c=1, \ldots, 10$,  the ${\Xb}_{2c;2j_2}$'s have concentration $5+c/10$ and still locations $\thetab_1=(1,0)\pr$ and $\thetab_2=(-1,0)\pr$.
  \item[(iii)] ${\Xb}_{3;1j_1}$ and ${\Xb}_{30;2j_2}$ have a common concentration $\kappa=10$ and locations $\thetab_1=(1,0)\pr$ and $\thetab_2=(-1,0)\pr$. Then for $c=1, \ldots, 10$,  the ${\Xb}_{3c;2j_2}$'s have concentration $10+c/10$ and still locations $\thetab_1=(1,0)\pr$ and $\thetab_2=(-1,0)\pr$.
\end{itemize}
For all $\ell=1,2,3$, the random vectors $\Xb_{\ell;1j_1}$ and $\Xb_{\ell 0;2j_2}$ are under the null hypothesis. Then, for $c=1, \ldots, 10$, the random vectors $\Xb_{\ell;1j_1}$ and $\Xb_{\ell 0;2j_2}$ are (increasingly) under the alternative. The results are plotted in Figures \ref{fighom1} and \ref{fighom2}. Inspection of the Figures reveals that the test $\phi_{\rm Hom}\n$ reaches the nominal level constraint even with small sample sizes. The power of the test decreases when the concentration increases.

\appendix

\section{Appendix: Proofs}

\noindent {\bf Proof of Proposition} \ref{alinea}
First note that, combining the fact that $\thetab$ and $\hat\thetab\n$ have norm 1 with the delta method applied to the mapping $\xb \mapsto \xb/  \| \xb \|$, we have that
\begin{eqnarray}\label{delta}
n^{1/2}(\hat{\thetab}- \thetab) &=& n^{1/2}\left(\frac{\hat{\thetab}}{\| \hat{\thetab}\|}- \frac{\thetab}{\| \thetab \|}\right) \nonumber \\
&=& ({\bf I}_k - \thetab \thetab\pr) n^{1/2}(\hat{\thetab}- \thetab)+ o_{\rm P}(1)
\end{eqnarray}
under ${\rm P}_\varthetab\n$ as $\ny$.
Now, the law of large numbers, the fact that ${\rm E}[\Xb_i]= {\rm E}(\Xb_i\pr \thetab) \thetab$ and (\ref{delta})  readily entail  that
\begin{eqnarray*}
T_{\hat \thetab}\n-T_\thetab\n &=& \bar{\Xb}\pr n^{1/2} (\hat{\thetab}- \thetab)\\
&=& ({\rm E}[\Xb_i])\pr ({\bf I}_k - \thetab \thetab\pr) n^{1/2}(\hat{\thetab}- \thetab)+ o_{\rm P}(1) \\
&=& {\rm E}(\Xb_i\pr \thetab) {\thetab}\pr ({\bf I}_k - \thetab \thetab\pr) n^{1/2}(\hat{\thetab}- \thetab)+ o_{\rm P}(1) \\
&=& o_{\rm P}(1)
\end{eqnarray*}
under ${\rm P}_\varthetab\n$ as $\ny$, which is the desired result.

\cqfd

\vspace{4mm}


\noindent {\bf Proof of Proposition} \ref{chideuxtwo} We first show Point (i). Since $\hat{\thetab}_i$ is a root-$n_i$ consistent estimator of $\thetab_i$, following the proof of Proposition \ref{alinea}, we readily have that $$ \bar{\Xb}_i\pr  n_i^{1/2} (\hat{\thetab}_i-\thetab_i)$$ is $o_{\rm P}(1)$ under ${\rm P}_{\varthetab^{(m)}}\n$ as $\ny$. As a direct consequence, letting $\hat{\varthetab}^{(m)}:=(\hat{\thetab}_1\pr, \ldots, \hat{\thetab}_m\pr, \hat{\kappa}, \ldots, \hat{\kappa})\pr$, we obviously have that
$$\Delta_{{\hat\varthetab}^{(m)}}^{({\rm I},i)(n)}-\Delta_{{\varthetab}^{(m)}}^{({\rm I},i)(n)}=-(r_i\n)^{1/2} n^{1/2} (A_k(\hat{\kappa})-A_k({\kappa}))+o_{\rm P}(1)$$
 under ${\rm P}_{\varthetab^m}\n$ with $\varthetab^m \in \mathcal{H}_0^{\rm Hom}$ as $\ny$. It follows directly that from the delta method and from Mardia and Jupp (2000, p. 199) that
 \begin{eqnarray} \label{alinbis} \Deltab_{{\hat\varthetab}^{(m)}}^{({\rm I})(n)}-\Deltab_{{\varthetab}^{(m)}}^{({\rm I})(n)} &=&- \Umb\n n^{1/2} (A_k(\hat{\kappa})-A_k({\kappa}))+o_{\rm P}(1)
 \\ &=&- \Umb\n A_k\pr(\kappa)\;  n^{1/2} (\hat{\kappa}-{\kappa})+o_{\rm P}(1) \nonumber \\
 &=&- \Umb\n D_k  \; n^{1/2} (\hat{\kappa}-{\kappa})+o_{\rm P}(1) \nonumber \\
 &=&- \Gamb_{\varthetab^{(m)}}^{({\rm II})} \Umb\n n^{1/2} (\hat{\kappa}-{\kappa})+o_{\rm P}(1) \nonumber
 \end{eqnarray}
 still under ${\rm P}_{\varthetab^{(m)}}\n$ as $\ny$ with $\varthetab^{(m)}$ such that all $\kappa_i$ components are equal (that is, we are under $\mathcal{H}_0^{\rm Hom}$). Therefore, defining
 $$(\Gamb_{\varthetab^{(m)}}^{({\rm I})})^{\perp}:=(\Gamb_{\varthetab^{(m)}}^{({\rm I})})^{-1}- \Umb\n \left((\Umb\n)\pr \Gamb_{\varthetab^{(m)}}^{({\rm I})}\Umb\n \right)^{-1} (\Umb\n)\pr,$$ 
 the consistency of $\hat{\varthetab}^{(m)}$ together with (\ref{alinbis}) entails that
 \begin{eqnarray*}
 Q_{\rm Hom}\n &=& Q_{\rm Hom}\n (\hat{\varthetab}^{(m)}) \\
 &=&  (\Deltab_{{\hat\varthetab}^{(m)}}^{({\rm I})})\pr (\Gamb_{\hat\varthetab^{(m)}}^{({\rm I})})^{\perp} \Deltab_{{\hat\varthetab}^{(m)}}^{({\rm I})} \\
 &=& (\Deltab_{{\varthetab}^{(m)}}^{({\rm I})}- \Gamb_{\varthetab^{(m)}}^{({\rm I})}\Umb\n n^{1/2} (\hat{\kappa}-{\kappa})))\pr (\Gamb_{\varthetab^{(m)}}^{({\rm I})})^{\perp} (\Deltab_{{\varthetab}^{(m)}}^{({\rm I})}- \Gamb_{\varthetab^{(m)}}^{({\rm I})}\Umb\n n^{1/2} (\hat{\kappa}-{\kappa})))+ o_{\rm P}(1) \nonumber \\
 &=& (\Deltab_{{\varthetab}^{(m)}}^{({\rm I})})\pr (\Gamb_{\varthetab^{(m)}}^{({\rm I})})^{\perp} \Deltab_{{\varthetab}^{(m)}}^{({\rm I})}+ o_{\rm P}(1) \\
 &=& Q_{\rm Hom}\n ({\varthetab}^{(m)})+o_{\rm P}(1)
 \end{eqnarray*}
 under ${\rm P}_{\varthetab^{(m)}}\n$ as $\ny$. Then, Point (i) directly follows from the asymptotic normality of $\Deltab_{{\varthetab}^{(m)}}^{({\rm I})}$ in Proposition \ref{ULANbis} and from the fact that $(\Gamb_{\varthetab^{(m)}}^{({\rm I})})^{\perp}\Gamb_{\varthetab^{(m)}}^{({\rm I})}$ is idempotent with trace $m-1$. Point (ii) follows by applying Le Cam's third Lemma as in Proposition \ref{chideuxthree}. For Point (iii), see Le Cam (1986) or Hallin and Paindaveine (2008) for a more recent reference. \cqfd

\vspace{1cm}

\noindent ACKNOWLEDGEMENTS

The research of Christophe Ley is supported by a Mandat de Charg\'e de Recherche from the Fonds
National de la Recherche Scientifique, Communaut\'e fran\c caise de Belgique.

\end{document}